\documentclass{amsart}

\input{epsf.tex}
\usepackage{amsfonts,amssymb,verbatim,amsmath,amsthm,latexsym,textcomp,amscd}
\usepackage{latexsym,amsfonts,amssymb,epsfig,verbatim, hyperref}
\usepackage{amsmath,amsthm,amssymb,latexsym,graphics,textcomp}
\usepackage{eucal}
\usepackage{graphicx}
\usepackage{color}
\usepackage{url}

\begin{document}

\newtheorem{theorem}{Theorem}[section]
\newtheorem{prop}[theorem]{Proposition}
\newtheorem{lemma}[theorem]{Lemma}
\newtheorem{cor}[theorem]{Corollary}
\newtheorem{definition}[theorem]{Definition}
\newtheorem{defn}[theorem]{Definition}
\newtheorem{conj}[theorem]{Conjecture}
\newtheorem{rmk}[theorem]{Remark}
\newtheorem{claim}[theorem]{Claim}
\newtheorem{defth}[theorem]{Definition-Theorem}
\newtheorem{qn}[theorem]{Question}

\newcommand{\boundary}{\partial}
\newcommand{\C}{{\mathbb C}}
\newcommand{\integers}{{\mathbb Z}}
\newcommand{\natls}{{\mathbb N}}
\newcommand{\ratls}{{\mathbb Q}}
\newcommand{\reals}{{\mathbb R}}
\newcommand{\proj}{{\mathbb P}}
\newcommand{\lhp}{{\mathbb L}}
\newcommand{\tube}{{\mathbb T}}
\newcommand{\cusp}{{\mathbb P}}
\newcommand\AAA{{\mathcal A}}
\newcommand\BB{{\mathcal B}}
\newcommand\CC{{\mathcal C}}
\newcommand\DD{{\mathcal D}}
\newcommand\EE{{\mathcal E}}
\newcommand\FF{{\mathcal F}}
\newcommand\GG{{\mathcal G}}
\newcommand\HH{{\mathcal H}}
\newcommand\II{{\mathcal I}}
\newcommand\JJ{{\mathcal J}}
\newcommand\KK{{\mathcal K}}
\newcommand\LL{{\mathcal L}}
\newcommand\MM{{\mathcal M}}
\newcommand\NN{{\mathcal N}}
\newcommand\OO{{\mathcal O}}
\newcommand\PP{{\mathcal P}}
\newcommand\QQ{{\mathcal Q}}
\newcommand\RR{{\mathcal R}}
\newcommand\SSS{{\mathcal S}}
\newcommand\TT{{\mathcal T}}
\newcommand\UU{{\mathcal U}}
\newcommand\VV{{\mathcal V}}
\newcommand\WW{{\mathcal W}}
\newcommand\XX{{\mathcal X}}
\newcommand\YY{{\mathcal Y}}
\newcommand\ZZ{{\mathcal Z}}
\newcommand\CH{{\CC\HH}}
\newcommand\MF{{\MM\FF}}
\newcommand\PMF{{\PP\kern-2pt\MM\FF}}
\newcommand\ML{{\MM\LL}}
\newcommand\PML{{\PP\kern-2pt\MM\LL}}
\newcommand\GL{{\GG\LL}}
\newcommand\Pol{{\mathcal P}}
\newcommand\half{{\textstyle{\frac12}}}
\newcommand\Half{{\frac12}}
\newcommand\Mod{\operatorname{Mod}}
\newcommand\Area{\operatorname{Area}}
\newcommand\ep{\epsilon}
\newcommand\hhat{\widehat}
\newcommand\Proj{{\mathbf P}}
\newcommand\U{{\mathbf U}}
 \newcommand\Hyp{{\mathbf H}}
\newcommand\D{{\mathbf D}}
\newcommand\Z{{\mathbb Z}}
\newcommand\R{{\mathbb R}}
\newcommand\Q{{\mathbb Q}}
\newcommand\E{{\mathbb E}}
\newcommand\til{\widetilde}
\newcommand\length{\operatorname{length}}
\newcommand\tr{\operatorname{tr}}
\newcommand\gesim{\succ}
\newcommand\lesim{\prec}
\newcommand\simle{\lesim}
\newcommand\simge{\gesim}
\newcommand{\simmult}{\asymp}
\newcommand{\simadd}{\mathrel{\overset{\text{\tiny $+$}}{\sim}}}
\newcommand{\ssm}{\setminus}
\newcommand{\diam}{\operatorname{diam}}
\newcommand{\pair}[1]{\langle #1\rangle}
\newcommand{\T}{{\mathbf T}}
\newcommand{\inj}{\operatorname{inj}}
\newcommand{\pleat}{\operatorname{\mathbf{pleat}}}
\newcommand{\short}{\operatorname{\mathbf{short}}}
\newcommand{\vertices}{\operatorname{vert}}
\newcommand{\collar}{\operatorname{\mathbf{collar}}}
\newcommand{\bcollar}{\operatorname{\overline{\mathbf{collar}}}}
\newcommand{\I}{{\mathbf I}}
\newcommand{\tprec}{\prec_t}
\newcommand{\fprec}{\prec_f}
\newcommand{\bprec}{\prec_b}
\newcommand{\pprec}{\prec_p}
\newcommand{\ppreceq}{\preceq_p}
\newcommand{\sprec}{\prec_s}
\newcommand{\cpreceq}{\preceq_c}
\newcommand{\cprec}{\prec_c}
\newcommand{\topprec}{\prec_{\rm top}}
\newcommand{\Topprec}{\prec_{\rm TOP}}
\newcommand{\fsub}{\mathrel{\scriptstyle\searrow}}
\newcommand{\bsub}{\mathrel{\scriptstyle\swarrow}}
\newcommand{\fsubd}{\mathrel{{\scriptstyle\searrow}\kern-1ex^d\kern0.5ex}}
\newcommand{\bsubd}{\mathrel{{\scriptstyle\swarrow}\kern-1.6ex^d\kern0.8ex}}
\newcommand{\fsubeq}{\mathrel{\raise-.7ex\hbox{$\overset{\searrow}{=}$}}}
\newcommand{\bsubeq}{\mathrel{\raise-.7ex\hbox{$\overset{\swarrow}{=}$}}}
\newcommand{\tw}{\operatorname{tw}}
\newcommand{\base}{\operatorname{base}}
\newcommand{\trans}{\operatorname{trans}}
\newcommand{\rest}{|_}
\newcommand{\bbar}{\overline}
\newcommand{\UML}{\operatorname{\UU\MM\LL}}
\newcommand{\EL}{\mathcal{EL}}
\newcommand{\tsum}{\sideset{}{'}\sum}
\newcommand{\tsh}[1]{\left\{\kern-.9ex\left\{#1\right\}\kern-.9ex\right\}}
\newcommand{\Tsh}[2]{\tsh{#2}_{#1}}
\newcommand{\qeq}{\mathrel{\approx}}
\newcommand{\Qeq}[1]{\mathrel{\approx_{#1}}}
\newcommand{\qle}{\lesssim}
\newcommand{\Qle}[1]{\mathrel{\lesssim_{#1}}}
\newcommand{\simp}{\operatorname{simp}}
\newcommand{\vsucc}{\operatorname{succ}}
\newcommand{\vpred}{\operatorname{pred}}
\newcommand\fhalf[1]{\overrightarrow {#1}}
\newcommand\bhalf[1]{\overleftarrow {#1}}
\newcommand\sleft{_{\text{left}}}
\newcommand\sright{_{\text{right}}}
\newcommand\sbtop{_{\text{top}}}
\newcommand\sbot{_{\text{bot}}}
\newcommand\sll{_{\mathbf l}}
\newcommand\srr{_{\mathbf r}}
\newcommand\geod{\operatorname{\mathbf g}}
\newcommand\mtorus[1]{\boundary U(#1)}
\newcommand\A{\mathbf A}
\newcommand\Aleft[1]{\A\sleft(#1)}
\newcommand\Aright[1]{\A\sright(#1)}
\newcommand\Atop[1]{\A\sbtop(#1)}
\newcommand\Abot[1]{\A\sbot(#1)}
\newcommand\boundvert{{\boundary_{||}}}
\newcommand\storus[1]{U(#1)}
\newcommand\Momega{\omega_M}
\newcommand\nomega{\omega_\nu}
\newcommand\twist{\operatorname{tw}}
\newcommand\modl{M_\nu}
\newcommand\MT{{\mathbb T}}
\newcommand\Teich{{\mathcal T}}
\renewcommand{\Re}{\operatorname{Re}}
\renewcommand{\Im}{\operatorname{Im}}

\title{Geometrically finite and infinite Kleinian Groups}

\author{Mahan Mj} 

\address{Indian Statistical Institute, 203, B.T. Road, Kolkata-700108}

\date{\today} 


\begin{abstract} This is a summary of the material for 3 lectures on geometrically finite and infinite Kleinian groups delivered by the author
at a conference held at Tata Institute of Fundamental Research in April 2014. 
\end{abstract}

\maketitle

\tableofcontents

\section{Lecture 1: Geometrically Finite Groups}
\subsection{Fuchsian  Groups}
A Kleinian group $G$ is a discrete subgroups of $PSL_2 ( {\mathbb{C}})
= Mob ({\widehat {\mathbb{C}}}) = Isom({\mathbb{H}}^3)$. This gives us three closely intertwined perspectives on the field:\\

\begin{enumerate}
\item Studying discrete subgroups $G$ of the group of Mobius transformations $Mob ({\widehat {\mathbb{C}}})$ emphasizes the
{\it Complex Analytical/Dynamic} aspect. 

\item 
Studying discrete subgroups $G$ of   $PSL_2 ( {\mathbb{C}})$ emphasizes the {\it Lie group/matrix group} theoretic aspect.

\item 
Studying discrete subgroups $G$ of   $Isom({\mathbb{H}}^3)$  emphasizes the {\it Hyperbolic Geometry} aspect.

\end{enumerate}

We shall largely emphasize the third perspective.
Since $G$ is discrete, we can pass to the quotient  $M^3 = {\mathbb{H}}^3/G$. Thus we are studying hyperbolic structures 
on $3$-manifolds.  

In order to obtain some examples, we first move one dimension down
and look at discrete subgroups $G$ of the group of Mobius transformations $Mob (\Delta ) = Mob ({\mathbb{H}})$ of the unit disk (which is
conformally equivalent to the upper half plane). 
These are called
Fuchsian Groups, and were discovered by Poincare.
The natural metric of constant negative curvature on the upper half
plane is given by  $ds^2 = \frac{dx^2 + dy^2}{y^2}$. This is called the {\it hyperbolic metric.} The resulting space is denoted
as ${\mathbb{H}}^2$.

The associated conformal structure is exactly the complex structure
on ${\mathbb{H}}  = \{ z \in \mathbb{C} : Im(z) > 0 \}$. It turns out that orientation preserving isometries of ${\mathbb{H}}^2$
are exactly the conformal automorphisms  of ${\mathbb{H}}^2$.
 The boundary circle $S^1$ compactifies $\Delta$. This has a geometric interpretation.
It codes the `ideal' boundary of ${\mathbb{H}}^2$, consisting of asymptote classes of geodesics.
The topology on $S^1$ is induced by a 
metric which is defined as the   angle subtended at $0 \in \Delta$.  The
geodesics turn out to be semicircles meeting the boundary $S^1$ at right angles.

We now proceed to construct an example of a discrete subgroup of $Isom({\mathbb{H}}^2)$. The genus two orientable surface
can be described as a quotient space of an octagon with  edges labelled $a_1, b_1, a_1^{-1}, b_1^{-1}, a_2, b_2,
a_2{-1}, b_2^{-1}$, where the boundary has the identification induced by this labelling. In order to construct a metric 
of constant negative curvature on it, we have to ensure that each point has a small neighborhood isometric to a small ball in 
${\mathbb{H}}^2$. To ensure this it is enough to do the above identification on  a regular hyperbolic octagon (all sides and all
angles equal) such that the sum of the interior angles is $2 \pi$.
To ensure this, we have to make each interior angle equal $\frac{2\pi}{8} $. The infinitesimal regular
octagon at the tangent space to the origin has interior angles equal to $\frac{3\pi}{4} $. Also the ideal 
regular
octagon in ${\mathbb{H}}^2$ has all interior angles zero. See figure below.

\begin{center}

\includegraphics[height=6cm]{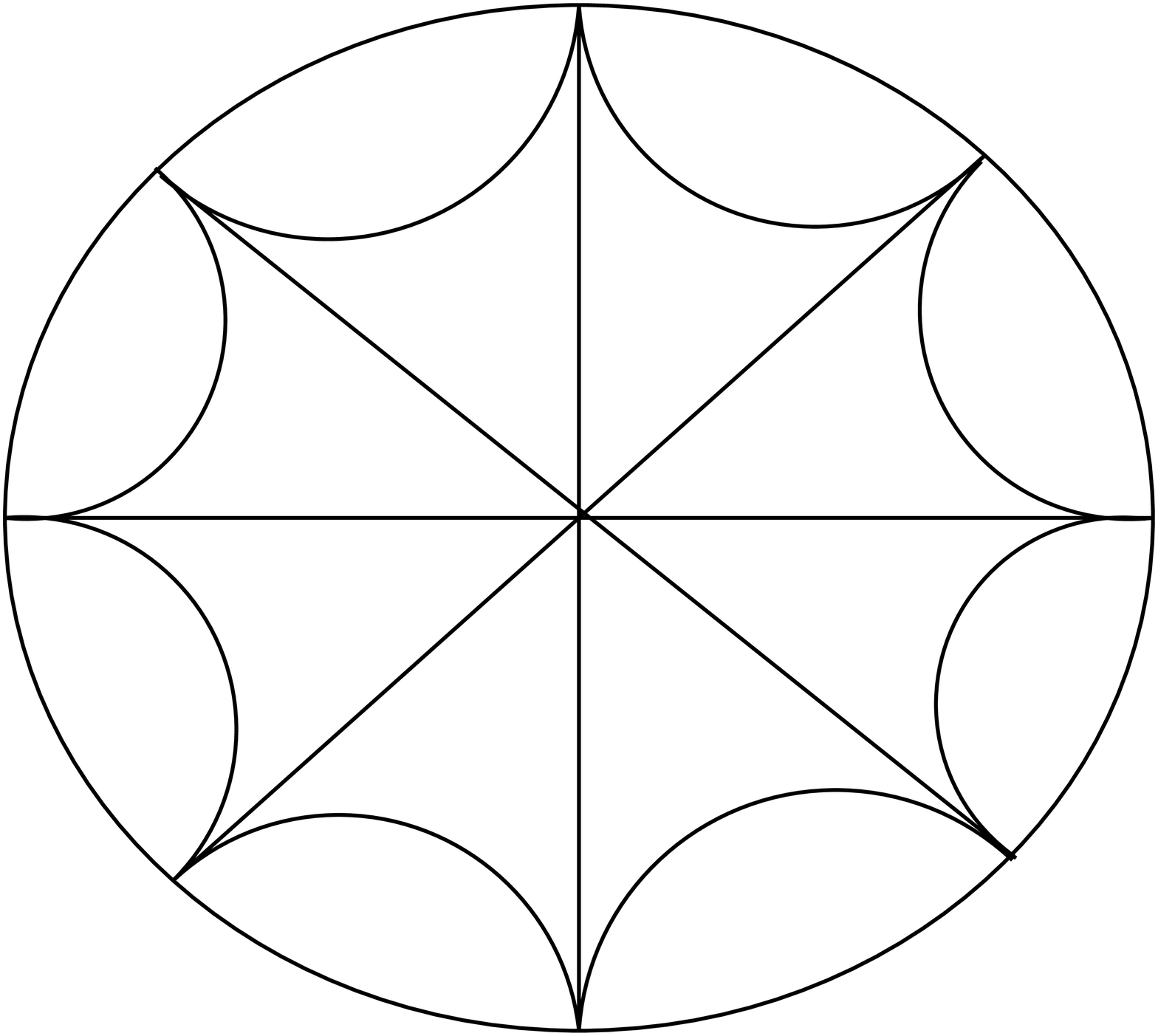}

\end{center}

Hence by the Intermediate value Theorem, as we increase the size of the octagon from an infinitesimal one to an ideal one,
we shall hit interior angles all equal to $\frac{\pi}{4} $ at some stage. The group $G$ that results from side-pairing transformations
corresponds to a Fuchsian group, or equivalently, a discrete faithful representation of the fundamental group
of a genus 2 surface into $Isom({\mathbb{H}}^2)$. We let $\rho$ denote the associated representation.

\subsection{Kleinian  Groups}

We now move back to ${\mathbb{H}}^3$. The hyperbolic metric is given by
 $ds^2 = \frac{dx^2 + dy^2 + dz^2}{z^2}$ on upper half space. Note that the
metric blows up as one approaches  $z=0$. Equivalently we could consider the ball model,
where the boundary $S^2 = {\widehat {\mathbb{C}}}$ consists of ideal end-points of geodesic rays as before. 
The metric on ${\widehat {\mathbb{C}}}$ is given by the angle subtended at $0 \in {\mathbb{H}}^3$.

Since $Isom({\mathbb{H}}^2) \subset Isom({\mathbb{H}}^3)$, we can look upon the discrete group $G$ we constructed above
also as a discrete subgroup of $Isom({\mathbb{H}}^3)$. 

\begin{center}

\includegraphics[height=6cm]{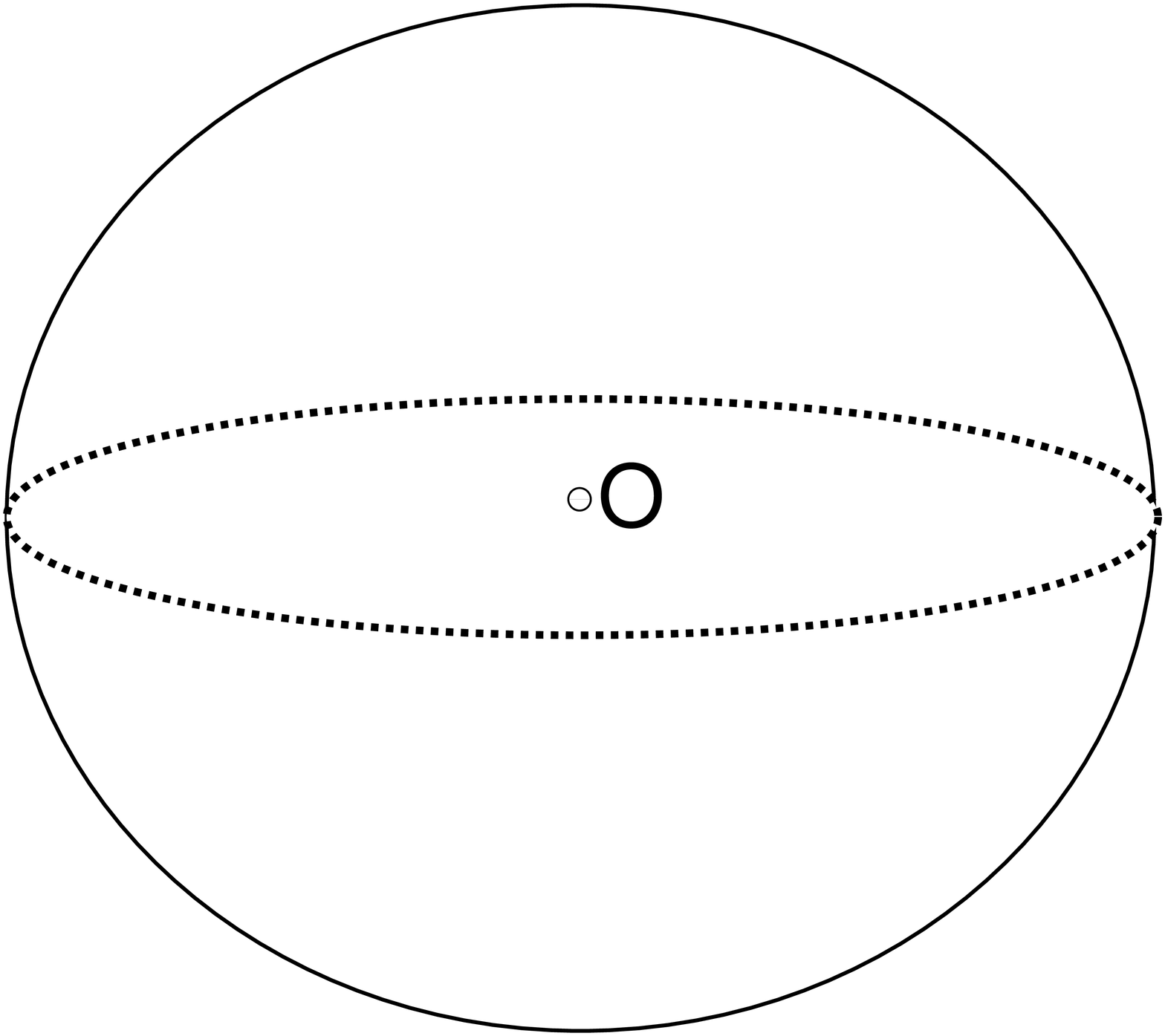}

\end{center}

In the above picture two things need to be observed. \\
1) the orbit $G.o$ accumulates on the equatorial circle. This is called the {\it limit set} $\Lambda_G$. \\
2) The complement of $\Lambda_G$ consists of two round open discs. 
On each of these disks, $G$ acts freely (i.e. without fixed points) properly discontinuously,
by conformal automorphisms. 
Hence quotient is two copies of the `same' Riemann surface (i.e. a one dimensional
complex analytic manifold). The complement
${\widehat {\mathbb{C}}} \setminus \Lambda_G = \Omega_G$ is called the {\it domain of discontinuity} of $G$.\\

We proceed with slightly more formal definitions identifying ${\widehat {\mathbb{C}}} $ with the sphere $S^2$.

\begin{defn} If $x \in \Hyp^3$ 
is any point, and $G$ is a discrete group of isometries, the limit set $\Lambda_G \subset S^2$
 is defined to be the set of accumulation points
of the orbit $G.x$ of $x$. 

The domain of discontinuity for a discrete group $G$ is defined to be $\Omega_G  = S^2 \setminus \Lambda_G $.\end{defn}

\begin{prop}\cite{thurstonnotes}[Proposition 8.1.2] If $G$ is not elementary, then every non-empty closed subset
of $S^2$ invariant by $G$ contains the limit set $\Lambda_G$. \end{prop}

Suppose that $G$ is abstractly isomorphic to the fundamental group of a finite area hyperbolic surface $S^h$,
and $\rho : \pi_1(S^h) \rightarrow PSL_2 ( {\mathbb{C}})$ be a representation with image $G$.
Suppose further that $\rho$ is {\it strictly type-preserving}, i.e. $g \in \pi_1(S^h)$ represents an element
in a peripheral (cusp) subgroup
  if and only $\rho (g)$ is parabolic. In this situation we shall refer to $G$ as a {\it surface Kleinian group}.
A recurring theme in the context of finitely generated, infinite covolume Kleinian groups is that the general
theory can be reduced to the study of surface Kleinian groups. Equivalently, we study the representation
space $Rep(\pi_1(S^h), PSL_2 ( {\mathbb{C}})$.

Regarding $G$ as a subgroup of $Mob({\widehat {\mathbb{C}}})$, the {\it dynamics} of the action
of $G$ on ${\widehat {\mathbb{C}}})$ emerges. The
{\it limit set} $\Lambda_G $ of $G$ is defined to be the set of accumulation points of the orbit
$G . o$ in ${\widehat {\mathbb{C}}}$
 for some (any) $o \in  {\mathbb{H}}^3$. The limit set is the locus of chaotic dynamics of the action
of $G$ on ${\mathbb{C}}$. The complement ${\widehat {\mathbb{C}}} \setminus \Lambda_G = \Omega_G$ is called the {\it domain of discontinuity} of $G$.

On the other hand regarding $G$ as a subgroup of $Isom({\mathbb{H}}^3)$, we obtain a quotient hyperbolic
3-manifold $M = {\mathbb{H}}^3/G$ with fundamental group $G$. 

A major problem in the theory of Kleinian groups is to {\bf understand the relationship between
the dynamic and the hyperbolic geometric descriptions} of $G$.

The { \bf Ahlfors-Bers simultaneous Uniformization Theorem} states that
given any two conformal structures $\tau_1, \tau_2$ on a surface, there is a discrete subgroup $G$ of 
$Mob ({\widehat {\mathbb{C}}})$ whose limit set is {\it topologically} a circle, and whose domain
of discontinuity quotients to the two Riemann surfaces $\tau_1, \tau_2$. See figure below.

\begin{center}

\includegraphics[height=6cm]{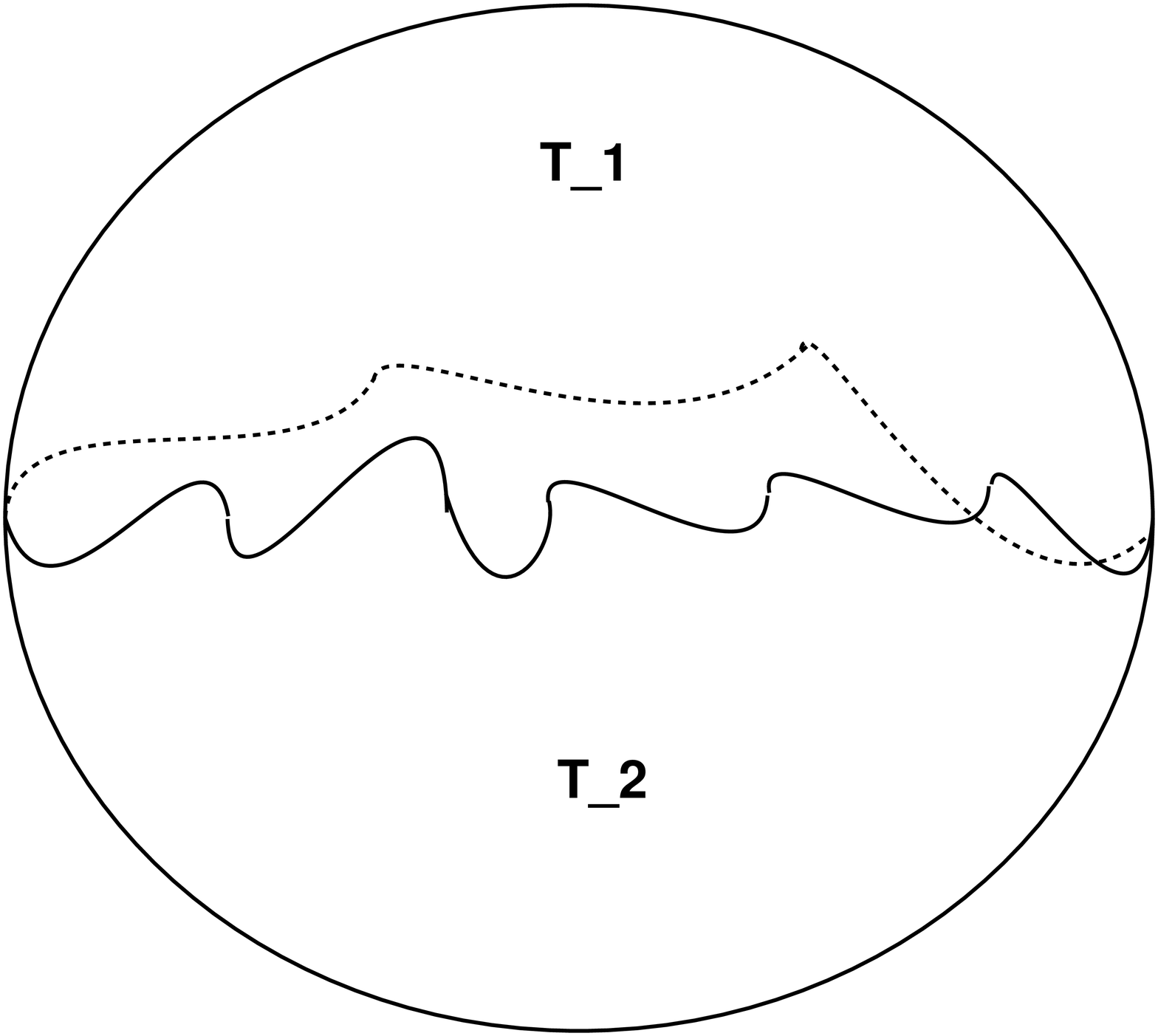}

\end{center}

The limit set is a quasiconformal
map of the round circle. These {\it (quasi Fuchsian)}
groups can be thought of as {\it deformations} of Fuchsian groups (Lie
group theoretically) or quasiconformal deformations (analytically). Ahlfors and Bers  proved that these are precisely
 all  quasiconvex surface Kleinian groups. 

The convex hull $CH_G$ of  $\Lambda_G$ is the smallest closed convex subset of ${\mathbb{H}}^3$ invariant under $G$. 
It can be constructed by joining all pairs of points on limit set by bi-infinite geodesics and iterating this
construction. The quotient of $CH_G$ by $G$, which is homeomorphic to $S^h \times [0,1]$, is
called the {\it Convex core $CC(M)$} of $M = {\mathbb{H}}^3/G$.

The `thickness' of $CC(M)$ for a quasi Fuchsian surface Kleinian group, 
measured by the distance between  $S^h \times \{ 0 \}$
and $S^h \times \{ 1 \}$ is a geometric measure of the complexity of the quasi Fuchsian group $G$.

\section{Lecture 2: Laminations} The main technical tools required to deal with the  notions of convex hulls and bending introduced In Lecture 1 are
laminations and pleated surfaces. I  followed Thurston's (unpublished) notes \cite{thurstonnotes} on the subject.

\begin{defn} A geodesic lamination on a hyperbolic surface is a foliation of a closed subset with geodesics. \end{defn}

Geodesic laminations arise naturally in a number of contexts in the study of hyperbolic 2- and 3- manifolds.

\begin{enumerate}
\item as stable and unstable laminations corresponding to a pseudo-anosov diffeomorphism of a hyperbolic surface.
\item as the pleating locus of a component of the convex core boundary $\partial CC(M)$ of a hyperolic 3-manifold $M$.
\item as the ending lamination corresponding to a geometrically infinite end of a hyperolic 3-manifold.
\end{enumerate}

We shall in this section discuss briefly how each of these examples arise.

\subsection{Stable and Unstable laminations}
We consider the torus $T^2$ equipped with a diffeomorphism $\phi$, whose action on homology is given by a $2 \times 2$ matrix with irrational
eigenvalues, e.g. $\frac{3 + \sqrt{5}}{2}, \frac{3 - \sqrt{5}}{2}$. 
Then the eigendirections give rise to two sets of foliations by dense copies of $\R$: the stable and unstable foliation. Such a diffeomorphism
is called Anosov. Anosov diffeomorphisms of the torus may be characterized in terms of their action on $\pi_1 (T)$ as not having periodic conjugacy classes.

Now consider the stable (or unstable) foliation (minus a point $\ast$) on $S = (T^2 \setminus \{ \ast \})$. Equip $S$ with a complete hyperbolic structure
of finite volume and straighten every leaf of the foliation to a complete geodesic. The resulting union of leaves is called the 
stable (or unstable) lamination of the diffeomorphism $\phi$ on the hyperbolic surface $S$.  

One of the fundamental pieces of Thurston's work \cite{FLP} shows that the existence of such a stable and unstable lamination generalizes to all hyperbolic
surfaces. A diffeomorphism $\phi$ of a hyperbolic surface $S$ preserving punctures (or boundary components according to taste)
is called pseudo Anosov if the action of $\phi_\ast$ on $\pi_1 (S)$ has no periodic conjugacy classes.
Thurston proved the existence of a unique stable and unstable lamination without any closed leaves for any pseudo Anosov 
diffeomorphism $\phi$ acting on a hyperbolic surface $S$.

\subsection{Pleating locus} We quote a picturesque passage from \cite{thurstonnotes}:

\begin{quote}
Consider a closed curve $\sigma$ in Euclidean space, and its convex hull $H(\sigma)$. The
boundary of a convex body always has non-negative Gaussian curvature. On the
other hand, each point $p$ in $\partial H(\sigma) \setminus \sigma$ lies in the interior of some line segment or
triangle with vertices on $\sigma$. Thus, there is some line segment on $\partial H(\sigma)$ through $p$,
so that $\partial H(\sigma)$ has non-positive curvature at $p$. It follows that $\partial H(\sigma) \setminus \sigma$ has zero
curvature, i.e., it is "developable". If you are not familiar with this idea, you can
see it by bending a curve out of a piece of stiff wire (like a coathanger). Now roll
the wire around on a big piece of paper, tracing out a curve where the wire touches.
Sometimes, the wire may touch at three or more points; this gives alternate ways
to roll, and you should carefully follow all of them. Cut out the region in the plane
bounded by this curve (piecing if necessary). By taping the paper together, you can
envelope the wire in a nice paper model of its convex hull. The physical process
of unrolling a developable surface onto the plane is the origin of the notion of the
developing map.

The same physical notion applies in hyperbolic three-space. If $K$ is any closed
set on $S^2$ (the sphere at infinity), then $H(K)$ is convex, yet each point on $\partial H(K)$ lies on a line segment
in $\partial H(K)$. Thus, $\partial H(K)$ can be developed to a hyperbolic plane. (In terms of
Riemannian geometry, $\partial H(K)$ has extrinsic curvature 0, so its intrinsic curvature
is the ambient sectional curvature, -1. Note however that $\partial H(K)$ is not usually
differentiable). Thus $\partial H(K)$ has the natural structure of a complete hyperbolic
surface.
\end{quote}

This forces $\partial H(K)$ equipped with its intrinsic metric to be a hyperbolic surface. However, there are complete geodesics
along which 
it is bent (but not crumpled). Thus each boundary component $S$, and hence its universal cover $\til S$, carries a metric that is intrinsically hyperbolic. However, in $\Hyp^3$,
the universal cover $\til S$ is bent along a geodesic lamination. $S$ is an example of a {\bf pleated surface}:

\begin{defn} \cite{thurstonnotes}[Definition 8.8.1] A {\bf pleated surface} in a hyperbolic three-manifold $N$ is
a complete hyperbolic surface $S$ of finite area, together with an isometric map $f :
S \to N$ such that every $x \in S$ is in the interior of some straight line segment which
is mapped by $f$ to a straight line segment. Also, $f$ must take every cusp of $S$ to a
cusp of $N$

The pleating locus of the pleated surface $f: S \to M$  is the set $\gamma \subset S$ consisting of those points in the pleated surface which are in the
interior of unique line segments mapped to line segments.
\end{defn}

\begin{prop} \cite{thurstonnotes}[Proposition 8.8.2] The pleating locus $\gamma$
 is a geodesic lamination on $S$. The map $f$ is totally geodesic in
the complement of 
$\gamma$. \end{prop}

\subsection{Ending Laminations} The notion of an ending lamination comes up in the context of a geometrically infinite group. We shall deal with these groups in greater detail in Lecture 3.  Thurston introduces the notion of a geometrically tame end $E$ of a manifold $M$. An end $E$ of a hyperbolic manifold
$M$ is geometrically tame (and geometrically infinite) if there exists a sequence of pleated surfaces exiting $E$.

For such an end $E$, choose a sequence of simple closed curves $\{ \sigma_n \}$ exiting $E$. Let $S = \partial E$ be the  bounding surface of $E$.
Then the limit of such a sequence (in a suitable sense; the reader
will not be much mistaken if (s)he thinks of the Hausdorff limit on the bounding surface $S$ of $E$) is a lamination $\lambda$. It turns out that
$\lambda$ is independent of the sequence $\{ \sigma_n \}$.

\section{Lecture 3: Geometrically Infinite Groups}
\subsection{Degenerate Groups}
The most intractable examples of surface Kleinian groups are obtained as limits of quasi Fuchsian groups.
In fact, it has been recently established by Minsky et al. \cite{minsky-elc1} \cite{minsky-elc2}  that
the set of all surface Kleinian groups (or equivalently all discrete faithful representations
of a surface group in $PSL_2({\mathbb{C}})$) are given by quasiFuchsian groups and their limits. This is known as the Bers density conjecture.

To construct limits of quasi Fuchsian groups, one allows the thickness of the convex core $CC(M)$ to tend to infinity. There are
two possibilities: \\
a) Let only $\tau_1$ degenerate. i.e. $I \rightarrow [0, \infty ) $ {\bf (simply degenerate case)}\\
b) Let both $\tau_1, \tau_2$ degenerate, i.e. $I \rightarrow (-\infty , \infty ) $ {\bf (doubly degenerate case)}\\

{\bf Thurston's Double Limit Theorem} \cite{thurston-hypstr2} says that these limits exist. A fundamental question in {\bf relating the
geometric and dynamic aspects} of Kleinian groups is the following.

\begin{qn} (Thurston) How does the limit set behave for the limiting manifold? \label{basic} \end{qn}

In the doubly degenerate case the limit set  is all of ${\widehat {\mathbb{C}}}$.

In the next section we outline our approach and solution to this problem.

\subsection{Extensions of Maps to Ideal Boundaries}
Starting with \cite{mitra-trees}, \cite{mitra-ct} and \cite{mitra-endlam}, we investigated the following question:

\begin{qn} Let $G$ be a hyperbolic group in the sense of Gromov acting freely and
properly discontinuously by isometries on a hyperbolic metric space
$X$. Does the inclusion of the Cayley graph $i: \Gamma_G \rightarrow
X$ extend continuously to the (Gromov) compactifications? \label{quest} \end{qn}

A positive answer to Question \ref{quest} gives us a precise handle on Question \ref{basic}.
In this generality the question first appears in \cite{mitra-thesis} (see also the Geometric Group Theory Problem
List \cite{bestvinahp}). As of date no counterexample is known.

However, special cases of Question \ref{quest} have been raised earlier in the context of Kleinian groups.

\noindent $\bullet 1$ In Section 6 of \cite{CT} (now
published as \cite{CTpub}), Cannon and Thurston
propose the following.

\begin{conj}
Suppose a  surface group $\pi_1 (S)$ acts freely and properly
discontinuously on ${\mathbb{H}}^3$ by isometries. Then the inclusion
$\tilde{i} : \widetilde{S^h} \rightarrow {\mathbb{H}}^3$ extends
continuously to the boundary \label{ct1} \end{conj}

The authors of \cite{CT} point out that for a simply degenerate group,
this is equivalent to asking if the limit set is locally connected.

\smallskip

\noindent $\bullet 2$ In \cite{ctm-locconn}, McMullen makes the
following more general conjecture: \\

\begin{conj}
For any hyperbolic 3-manifold $N$ with finitely generated fundamental
group, there exists a continuous, $\pi_1(N)$-equivariant map \\
\begin{center}
$F: \partial \pi_1 (N) \rightarrow \Lambda \subset S^2_{\infty}$

\end{center}

where the boundary $\partial \pi_1(N)$ is constructed by scaling the
metric on the Cayley graph of $\pi_1 (N)$ by the conformal factor of
$d(e,x)^{-2}$, then taking the metric completion. (cf. Floyd
\cite{Floyd})
\label{ct2} \end{conj}

\smallskip

In \cite{mahan-split} and \cite{mahan-kl} we provide a complete positive answer to both Conjectures \ref{ct1} and \ref{ct2}.

As a consequence we also establish in \cite{mahan-split} the following Theorem which proves a long-standing conjecture in the theory of Kleinian groups \cite{abikoff}
\cite{CT}.

\begin{theorem} Connected limit sets of finitely generated Kleinian groups is locally connected. \label{lc} \end{theorem}

In the next subsection, after describing the history of these problems, we shall give more details about the structure of limit sets and their relation to the
geometry of surface Kleinian groups.

\subsection{History and Solution of the Problem}

In \cite{abikoff}, Abikoff (1976) claimed to prove that limit sets of simply degenerate surface Kleinian groups were never locally connected. Thurston and Kerckhoff 
found a flaw in his proof in about 1980.

The first major result that started this entire program was Cannon and
Thurston's result  \cite{CT} for hyperbolic 3-manifolds fibering over the
circle with fiber a closed surface group. 

Let $M$ be a closed hyperbolic 3-manifold fibering over the circle with 
fiber $F$. Let $\widetilde F$ and $\widetilde M$ denote the universal
covers of $F$ and $M$ respectively. Then $\widetilde F$ and $\widetilde M$
are quasi-isometric to ${\mathbb{H}}^2$ and ${\mathbb{H}}^3$ respectively. Now let
${{\mathbb{D}}^2}={\mathbb{H}}^2\cup{\mathbb{S}}^1_\infty$ and 
${{\mathbb{D}}^3}={\mathbb{H}}^3\cup{\mathbb{S}}^2_\infty$
denote the standard compactifications. In \cite{CT} Cannon and Thurston
show that the usual inclusion of $\widetilde F$ into $\widetilde M$
extends to a continuous map from ${\mathbb{D}}^2$ to ${\mathbb{D}}^3$.
This was extended to Kleinian surface groups of bounded geometry
without parabolics by Minsky \cite{minsky-jams}.

An
alternate approach (purely in terms of coarse geometry ignoring all
local information) was given by the author in \cite{mitra-trees}
generalizing the results of both  Cannon-Thurston and Minsky. We proved the
Cannon-Thurston result for hyperbolic 3-manifolds of bounded geometry
without parabolics and with freely indecomposable fundamental group. A
different approach based on Minsky's work was
given by Klarreich \cite{klarreich}. 

Bowditch \cite{bowditch-ct}
\cite{bowditch-stacks} proved the Cannon-Thurston result
for punctured surface Kleinian groups of
bounded geometry.
In  \cite{brahma-pared} we gave an
alternate proof of Bowditch's results and simultaneously generalized the
results of Cannon-Thurston, Minsky, Bowditch, and those of 
\cite{mitra-trees} to
all 3 manifolds of bounded geometry whose cores are incompressible
away from cusps.  The proof has the
advantage that it reduces to a proof for manifolds without parabolics
when the 3 manifold in question has freely indecomposable fundamental
group and no accidental parabolics.

In the expository paper \cite{brahma-bddgeo} we give our proof of the results of Cannon and Thurston \cite{CT}, Minsky \cite{minsky-jams},
and Bowditch \cite{bowditch-ct} using the ideas of \cite{mitra-trees} and \cite{brahma-pared}.

In \cite{minsky-torus} Minsky established a bi-Lipschitz model for all punctured torus Kleinian groups. 
McMullen \cite{ctm-locconn} proved the Cannon-Thurston result
for punctured torus groups, using
 Minsky's model for these groups \cite{minsky-torus}. 

 In \cite{mahan-ibdd} we identified a large-scale 
coarse geometric structure involved in the Minsky model for punctured
torus groups (and called it {\bf
  i-bounded geometry}). {\em i-bounded geometry} can roughly be
regarded as that geometry of ends where the boundary tori of Margulis
tubes have uniformly bounded diameter. We gave a  proof for models of
{\it i-bounded geometry}. In combination with the methods of
\cite{brahma-pared} this was enough to bring under the same umbrella
all known results on Cannon-Thurston maps for 3 manifolds whose cores are
incompressible away from cusps. 

In \cite{brahma-amalgeo} we further generalized possible geometries allowing us to push our techniques through to establish
the Cannon-Thurston property.

In the mean time, in the proof of the celebrated Ending Lamination Conjecture, Minsky \cite{minsky-elc1} and Brock-Canary-Minsky \cite{minsky-elc2}
established a bi-Lipschitz model for all surface Kleinian groups.

In \cite{mahan-split}, we used the Minsky model of \cite{minsky-elc1} to prove that all hyperbolic
3-manifolds homotopy equivalent to a surface satisfy the conditions imposed in the geometries dealt with in \cite{brahma-amalgeo}.
This establishes the Cannon-Thurston property for all surface Kleinian groups and proves Conjecture \ref{ct1}.
It follows that surface Kleinian groups have locally connected limit sets. Combining this result with a reduction Theorem of Anderson and Maskit \cite{and-mask},
we prove that connected limit sets of  finitely generated Kleinian groups are locally connected (Theorem \ref{lc}).
Finally in \cite{mahan-kl} we extend the techniques of \cite{mahan-split} to cover handlebody groups and prove Conjecture \ref{ct2}.

We then gave explicit descriptions of the boundary identifications of \cite{mahan-split} in terms of {\it ending laminations} in \cite{mahan-elct}. This finally yields a rather complete and satisfactory solution to Question \ref{basic}.

\bibliography{tifr2015}
\bibliographystyle{alpha}

\end{document}